\documentclass[12pt]{amsart}
\newtheorem{Theorem}{Theorem}[section]

\newtheorem{Lemma}[Theorem]{Lemma}

\theoremstyle{definition}

\theoremstyle{remark}

\begin{document}
\sloppy
\title{A class of reconstructible graphs}
\author{Tetsuya Hosaka} 
\address{Department of Mathematics, Utsunomiya University, 
Utsunomiya, 321-8505, Japan}
\date{July 11, 2004}
\email{hosaka@cc.utsunomiya-u.ac.jp}
\subjclass[2000]{}
\thanks{Partly supported by the Grant-in-Aid for Scientific Research, 
The Ministry of Education, Culture, Sports, Science and Technology, Japan, 
(No.~15740029).}
\maketitle
\begin{abstract}
In this paper, we give a class of reconstructible graphs.
\end{abstract}

\section{Introduction}

The purpose of this paper is to give a class of reconstructible graphs.

For a graph $G$ and $v\in V(G)$, 
let $N(v)=\{u\in V(G)\,|\, uv\in E(G)\}$ and 
$N[v]=N(v)\cup\{v\}$.

We prove the following theorem.

\begin{Theorem}
Let $G$ be a graph with $V(G)=\{v_1,\dots,v_k,\dots,v_n\}$ $(1<k<n)$.
Suppose that 
\begin{enumerate}
\item[${\rm (i)}$] $\bigcup_{i=1}^{k}N[v_i]=V(G)$,
\item[${\rm (ii)}$] $|d(v_i)-d(v_j)|>1$ for any $i\in\{1,\dots,k\}$ 
and $j\in\{1,\dots,n\}\setminus\{i\}$ and
\item[${\rm (iii)}$] $N[v_1]\cap \bigcup_{i=2}^{k}N[v_i]=\emptyset$.
\end{enumerate}
Then $G$ is reconstructible.
\end{Theorem}

\section{Proof of the theorem}

Let $G$ be a graph with $V(G)=\{v_1,\dots,v_k,\dots,v_n\}$ $(1<k<n)$.
Suppose that 
\begin{enumerate}
\item[${\rm (i)}$] $\bigcup_{i=1}^{k}N[v_i]=V(G)$,
\item[${\rm (ii)}$] $|d(v_i)-d(v_j)|>1$ for any $i\in\{1,\dots,k\}$ 
and $j\in\{1,\dots,n\}\setminus\{i\}$ and
\item[${\rm (iii)}$] $N[v_1]\cap \bigcup_{i=2}^{k}N[v_i]=\emptyset$.
\end{enumerate}
Let $G'$ be a graph with $V(G')=\{v'_1,\dots,v'_n\}$ 
such that $G-v_j \cong G'-v'_j$ for any $j\in\{1,\dots,n\}$.

Here we note that $d(v_j)=d(v'_j)$ for any $j\in\{1,\dots,n\}$.

Let $f_j:G-v_j \rightarrow G'-v'_j$ be an isomorphism 
for each $j\in\{1,\dots,n\}$.

\begin{Lemma}\label{lem1}
For each $i\in\{1,\dots,k\}$ and $j\in\{1,\dots,n\}\setminus\{i\}$, 
\begin{enumerate}
\item[(1)] $f_j(v_i)=v'_i$ and 
\item[(2)] $v_iv_j\in E(G)$ if and only if $v'_iv'_j\in E(G')$.
\end{enumerate}
\end{Lemma}

\begin{proof}
(1) We obtain that $f_j(v_i)=v'_i$ from ${\rm (ii)}$.

(2) Suppose that $v_iv_j\in E(G)$.
Then $d(v_i;G-v_j)=d(v_i)-1$, 
where $d(v_i;G-v_j)$ is the degree of $v_i$ in $G-v_j$.
Here either $d(v_i';G'-v'_j)=d(v'_i)$ or $d(v_i';G'-v'_j)=d(v'_i)-1$.
Since $G-v_j\cong G'-v'_j$, 
$d(v_i';G'-v'_j)=d(v'_i)-1$ by (ii).
This means that $v'_iv'_j\in E(G')$.
Thus if $v_iv_j\in E(G)$ then $v'_iv'_j\in E(G')$. 
By the same argument, we also can show the reverse.
\end{proof}

\begin{Lemma}\label{lem2}
~
\begin{enumerate}
\item[${\rm (i')}$] $\bigcup_{i=1}^{k}N[v'_i]=V(G')$.
\item[${\rm (ii')}$] $|d(v'_i)-d(v'_j)|>1$ for any $i\in\{1,\dots,k\}$ 
and $j\in\{1,\dots,n\}\setminus\{i\}$.
\item[${\rm (iii')}$] $N[v'_1]\cap \bigcup_{i=2}^{k}N[v'_i]=\emptyset$.
\end{enumerate}
\end{Lemma}

\begin{proof}
${\rm (ii')}$ is obvious because 
$d(v_j)=d(v'_j)$ for any $j\in\{1,\dots,n\}$.

${\rm (i')}$ Let $j\in\{1,\dots,n\}$. 
If $j\le k$ then $v'_j\in N[v'_j]$ and 
$v'_j\in \bigcup_{i=1}^{k}N[v'_i]$.
Suppose that $k+1\le j$. 
Since $v_j\in V(G)=\bigcup_{i=1}^{k}N[v_i]$ by (i), 
there exists $i_0\in\{1,\dots,k\}$ such that $v_j\in N[v_{i_0}]$.
Then $v_{i_0}v_j\in E(G)$.
By Lemma~\ref{lem1}~(2), 
$v'_{i_0}v'_j\in E(G')$, i.e., 
$v'_j\in N[v'_{i_0}]\subset \bigcup_{i=1}^{k}N[v'_i]$.
Thus $\bigcup_{i=1}^{k}N[v'_i]=V(G')$.

${\rm (iii')}$ 
Suppose that 
there exists $j\in\{1,\dots,n\}$ 
such that $v'_j\in N[v'_1]\cap \bigcup_{i=2}^{k}N[v'_i]$.
Then $v'_j\in N[v'_1]\cap N[v'_{i_0}]$ for some $i_0\in\{2,\dots,k\}$.
This means that $v'_1v'_j\in E(G')$ and $v'_{i_0}v'_j\in E(G')$.
By Lemma~\ref{lem1}~(2), 
$v_1v_j\in E(G)$ and $v_{i_0}v_j\in E(G)$.
Hence $v_j\in N[v_1]\cap N[v_{i_0}]$.
This contradicts ${\rm (iii)}$.
Thus $N[v'_1]\cap \bigcup_{i=2}^{k}N[v'_i]=\emptyset$.
\end{proof}

Now $f_1:G-v_1\rightarrow G'-v'_1$ is an isomorphism.
We show that the map $f:G\rightarrow G'$ defined by 
$f(v_1)=v'_1$ and $f(v_j)=f_1(v_j)$ for each $j\in\{2,\dots,n\}$ 
is an isomorphism.
It is sufficient to prove that 
for each $j\in\{2,\dots,n\}$, 
$v_1v_j\in E(G)$ if and only if $v'_1f_1(v_j)\in E(G')$.
Here we note that $v_1v_j\in E(G)$ if and only if $v_j\in N(v_1)$, and 
$v'_1f_1(v_j)\in E(G')$ if and only if $f_1(v_j)\in N(v'_1)$.
Hence we show that $f_1(N(v_1))=N(v'_1)$.
For each $i\in\{2,\dots,k\}$, by ${\rm (iii)}$ and ${\rm (iii')}$, 
\begin{align*}
f_1(N[v_i])&=f_1(N[v_i;G-v_1])=N[f_1(v_i);G'-v'_1] \\
&=N[v'_i;G'-v'_1]=N[v'_i].
\end{align*}
Hence 
\begin{align*}
f_1(N(v_1))&=f_1(V(G-v_1)-\bigcup_{i=2}^{k}N[v_i]) \\
&=V(G'-v'_1)-\bigcup_{i=2}^{k}f_1(N[v_i]) \\
&=V(G'-v'_1)-\bigcup_{i=2}^{k}N[v'_i] \\
&=N(v'_1) 
\end{align*}
by ${\rm (i)}$, ${\rm (iii)}$, ${\rm (i')}$ and ${\rm (iii')}$.
Thus $f_1(N(v_1))=N(v'_1)$ and 
$f:G\rightarrow G'$ is an isomorphism.
Therefore $G$ is reconstructible.

%

%
\end{document}